\documentclass[a4paper]{amsart}

\usepackage{amsmath,amssymb,amsthm,epsfig}

\def\dis
{\displaystyle}
\def\e
{\varepsilon}

\def\R{{\mathbb R}}
\def\N{{\mathbb N}}

\def\F{{\mathcal F}}
\def\S{{\mathcal S}}
\def\virgp{\raise 2pt\hbox{,}}
\def\bv{{\bf v}}
\def\bu{{\bf u}}
\def\bw{{\bf w}}
\def\1{{\rm 1\mskip-4.5mu l} }

\def\<{\langle}
\def\>{\rangle}

\def\({\left(}
\def\){\right)}

\def\O{\mathcal{O}}

\def\Eq#1#2{\mathop{\sim}\limits_{#1\rightarrow#2}}
\def\Tend#1#2{\mathop{\longrightarrow}\limits_{#1\rightarrow#2}}

\def\d{{\partial}}
\def\g{{\gamma}}
\def\h{{\hbar}}
\def\e{\varepsilon}
\def\si{{\sigma}}

\theoremstyle{plain}
\newtheorem{theo}{Theorem}[section]

\newtheorem{prop}[theo]{Proposition}
\newtheorem{hyp}[theo]{Assumptions}

\theoremstyle{definition}
\newtheorem{defin}[theo]{Definition}

\theoremstyle{remark}
\newtheorem{rema}[theo]{Remark}

\numberwithin{equation}{section}

\begin{document}
\title[Cascade of phase shifts for NLS]{Cascade of phase shifts
for nonlinear Schr\"odinger equations}
\author[R. Carles]{R{\'e}mi Carles}
\address{MAB, UMR CNRS 5466\\
Universit{\'e} Bordeaux 1\\ 351 cours de la Lib{\'e}ration\\ 33~405 Talence
cedex\\ France}
\email{Remi.Carles@math.cnrs.fr}
\begin{abstract}
We consider a semi-classical nonlinear Schr\"odinger equation. For
initial data causing focusing at one point in the linear case, we
study a nonlinearity which is super-critical in terms of asymptotic
effects near the caustic. We prove the existence of infinitely many
phase shifts appearing at the approach of the critical time. This phenomenon
is suggested by a formal computation. The rigorous proof shows a
quantitatively different asymptotic behavior. 
We explain these aspects, and discuss some problems left open.  
\end{abstract}
\thanks{This work was initiated at CMAF
(University of Lisbon), 
  partly supported by the FCT grant SFRH/BPD/16766/2004, and completed
  at IRMAR (University of Rennes). The author would like to thank these
  institutions for their hospitality. Support
  by the European network HYKE,  
funded  by the EC as contract  HPRN-CT-2002-00282, is also acknowledged.}
\maketitle

\section{Introduction}
\label{sec:intro}
We consider the semi-classical limit of the Cauchy problem, for
$(t,x)\in \R_+\times \R^n $:
\begin{equation}\label{eq:1}
  i\e \d_t \bu^\e + \frac{\e^2}{2} \Delta \bu^\e = f\(
  |\bu^\e|^{2}\)\bu^\e \quad ; \quad 
  \bu^\e_{\mid t=0}= \e^{k/2}a_0(x) e^{-i\frac{|x|^2}{2\e}} .  
\end{equation}
In the linear case $f\equiv 0$, the quadratic oscillations of the
initial data cause focusing at the origin at time $t=1$ in the limit
$\e \to 0$ (see Section~\ref{sec:free}). 
In the nonlinear case, the effective nonlinear effects
strongly depend on the size of the initial data, that is on
$k$. Changing notations, we consider:
\begin{equation}
  \label{eq:princ}
  i\e \d_t u^\e + \frac{\e^2}{2} \Delta u^\e = f\(\e^k
  |u^\e|^{2}\)u^\e \quad ; \quad 
  u^\e_{\mid t=0}= a_0(x) e^{-i\frac{|x|^2}{2\e}} .  
\end{equation}
In \cite{Ca2}, we justified  the general
heuristics presented in \cite{HK87}, in the case of \eqref{eq:princ},
for $f$ homogeneous of degree $\si$, $f(y)=y^\si$. Two notions of
criticality exist for $k$: outside the focal point, and near the focal
point, where the amplitude of $u^\e$ is strongly modified. We
described the sub-critical and critical cases. The aim of the present
paper is to study a supercritical case. 

Consider the case $f(y)=y^\si$, and denote $\alpha = k\si$. If $a_0\in
H^1(\R^n)$ with $|x|a_0 \in L^2(\R^n)$ and 
$\si<2/(n-2)$ when $n\ge 3$, $u^\e$ is defined globally in time in
$H^1(\R^n)$. The
following distinctions were established in \cite{Ca2}:
\begin{center}
\begin{tabular}[c]{l|c|c}
 &$\alpha >n\si $  & $\alpha=n\si$ \\
\hline
$\alpha >1$ &linear caustic&nonlinear caustic
\\
&linear WKB&linear WKB\\
\hline
$\alpha=1 $&linear caustic&nonlinear caustic\\
&nonlinear WKB&nonlinear WKB\\
\end{tabular}
\end{center}

\medskip

The term ``linear WKB'' means that outside the caustic, the
propagation of $u^\e$ can be described by a geometrical optics
approximation, with only linear effects involved at leading order. The
term ``linear caustic'' means that nonlinear effects are negligible at
leading order when the solution crosses the focal point. In either of
the two critical cases, nonlinear phenomena are described (the doubly
critical case was studied more precisely in \cite{CaCMP}); we recall
the case ``nonlinear caustic, linear WKB'' in
Section~\ref{sec:critical}. The case we study in this 
paper corresponds to $n\si >\alpha> 1$: super-critical caustic with linear
WKB r\'egime.

A look at conservation laws suggests the existence of new relevant
scales. In the case discussed so far, the conservations of charge and
energy read (see e.g. \cite{CazCourant}):
\begin{equation}
  \label{eq:conservation}
  \begin{aligned}
   & \|u^\e(t)\|_{L^2} = \|a_0\|_{L^2}\, ,\\
\frac{1}{2}&\|\e\nabla_x u^\e(t)\|_{L^2}^2 + \frac{\e^\alpha}{\si +1}
\| u^\e(t)\|_{L^{2\si +2}}^{2\si +2} = \text{const.}=\O (1) \Eq \e 0
\frac{1}{2}\|x a_0\|_{L^2}^2\, .   
  \end{aligned}
\end{equation}
When $\alpha \ge n\si$, the boundedness of $u^\e$ and $\e\nabla_x
u^\e$ in $L^2$ implies, along with Gagliardo--Nirenberg inequalities,
for $0<\e \le 1$:
\begin{equation*}
  \e^\alpha
\| u^\e(t)\|_{L^{2\si +2}}^{2\si +2}\le \e^{n\si}\| u^\e(t)\|_{L^{2\si
    +2}}^{2\si +2} \lesssim \|u^\e(t)\|_{L^2}^{2
  -(n-2)\si}\|\e \nabla_x u^\e(t)\|_{L^2}^{n\si}=\O(1)\,. 
\end{equation*}
Thus, linear
arguments allow us to recover a control of the nonlinear term in the
energy. Such a line of reasoning fails when $\alpha <n\si$: the control
provided by the conservation of energy hides stronger nonlinear
effects. In the linear case $f\equiv 0$, and $a_0\in \S(\R^n)$, one
can check that the following point-wise estimate holds:
\begin{equation}\label{eq:estlinsharp}
  |u^\e(t,x)|\lesssim \frac{1}{\(\e +|t-1|\)^{n/2}}\, \cdot
\end{equation}
In the four cases of the table, the same
estimate holds for the nonlinear solution in space dimension one
(\cite{Ca2,CaCMP}); like in the linear case, it is sharp. Thus, the
above estimate $\e^{n\si}\| 
u^\e(t)\|_{L^{2\si 
    +2}}^{2\si +2} =\O(1)$ is sharp only near the focal point. We now use
the \emph{a priori} estimate $\e^\alpha\| u^\e(t)\|_{L^{2\si
    +2}}^{2\si +2} =\O(1)$ given by the conservation of energy only
for $t\thickapprox 1$. Assuming that like in all the cases of the above
table, at time $t=1$, $u^\e$ is described by a concentrating
profile, 
\begin{equation*}
  u^\e(t,x) \Eq \e 0 \frac{1}{\e^{n\g/2}}\phi\( \frac{x}{\e^\g}\)\, ,
\end{equation*}
we check that the ``linear'' value $\g =1$ is forbidden (the power of
$\e$ in front of $\phi$ is to ensure the $L^2$-norm
conservation). Guessing that the nonlinear term $\e^\alpha\| u^\e(t)\|_{L^{2\si
    +2}}^{2\si +2}$ in the energy is exactly of order $\O(1)$ at the
caustic, we find $\g =\alpha/(n\si)$, that is:
\begin{equation}
  \label{eq:gamma}
  \g = \frac{k}{n}<1\, .
\end{equation}
We will not prove that the above argument is correct (see
Section~\ref{sec:beyond}), but we will show that the scale $\e^\g$ is an
important feature of this problem. Notice also that the above argument
suggests that the amplification of the solution $u^\e$ as time goes to
$1$ is less important than in the linear case; super-critical
phenomena may occur in the phase, and also affect the amplitude.  

We now go back to the notation \eqref{eq:princ}, and do not assume in
general that the nonlinearity is homogeneous (unless it is cubic):
\begin{hyp}\label{hyp:tout}
  The space dimension is $n\ge 2$.\\
  The initial amplitude belongs to the Schwartz space: $a_0\in \S(\R^n)$.\\
  The nonlinearity is smooth: $f\in
  C^\infty(\R_+;\R)$. \\
  $f(0)=0$ and $f'>0$. In particular, the
  nonlinearity is cubic at the origin. 
\end{hyp}
\begin{rema}\label{rema:reg}
  We suppose $a_0\in \S(\R^n)$ to avoid to count derivatives when not
  necessary. We could as well assume that $a_0$ belongs to Sobolev
  type spaces. If we require a control on the growth
  of $f$ at infinity, 
$(0\le) f(y)\lesssim \< y\>^q$ for $q<\frac{2}{n-2}$ when $n\ge 3$, 
then for every fixed $\e>0$, $u^\e$ is global in time, continuous with
  values in $H^1(\R^n)$ (see e.g. \cite{CazCourant}). This includes
  a cubic nonlinearity   in space dimension two or three.
\end{rema}
\begin{rema}
The assumption $f(0)=0$ is only to simplify notations,
since replacing $f$ with $f-f(0)$ turns $u^\e(t,x)$ into
$u^\e(t,x)e^{if(0)t/\e}$.  
\end{rema}
\begin{rema}
The assumption of the nonlinearity being cubic at the origin is
reminiscent of the paper by E.~Grenier
\cite{Grenier98} (see also P.~G\'erard \cite{PGX93}).  The proof of
our main result  
relies on ideas introduced in \cite{Grenier98} (see
Section~\ref{sec:justif}). 
\end{rema}  
\begin{rema}
  The one-dimensional cubic nonlinear Schr\"odinger equation is
  integrable. The 
  case $k=0$ with more general WKB data was treated in \cite{JLM}. 
\end{rema}
Before stating our main result, we give the following definition (see
e.g. \cite{RauchUtah}):
\begin{defin}
  If $T>0$, $(k_j)_{j\ge 1}$ is an increasing sequence of real
  numbers, $(\phi_j)_{j\ge 1}$ is a sequence in
  $H^{\infty}(\R^n):=\cap_{s\ge 0}H^s(\R^n)$, and $\phi
  \in C([0,T];H^s(\R^n))$ for every $s>0$, the asymptotic 
  relation  
  \begin{equation*}
    \phi(t,x)\sim \sum_{j\ge 1}t^{k_J}\phi_j(x)\quad \text{as }t\to 0
  \end{equation*}
means that for every integer $J\ge 1$ and every $s>0$, 
\begin{equation*}
  \left\| \phi(t,\cdot)-\sum_{j= 1}^J t^{k_j}\phi_j
  \right\|_{H^s(\R^n)} =o\(t^{k_J}\)\quad \text{as }t\to 0\, .
\end{equation*}
\end{defin}
\begin{theo}\label{theo:main}
Let Assumptions~\ref{hyp:tout} be satisfied. Assume $n>k>1$. Then
 there exist $T>0$ independent of $\e\in ]0,1]$, 
 a sequence $(\phi_j)_{j\ge 1}$ in $H^{\infty}(\R^n)$, and $\phi \in
C([0,T];H^s(\R^n))$ for every $s>0$, such that:\\
$1.$ $\phi(t,x)\sim \sum_{j\ge 1}t^{jn-1}\phi_j(x)$ as $t\to 0$.\\
$2.$ For $1-t \gg \e^\g$ ($\g=k/n<1$), the
asymptotic behavior of $u^\e$ is given by:
\begin{align*}
 &\limsup_{\e \to 0}\sup_{0\le t\le 1-\Lambda \e^\g}\left\| u^\e(t)-v^\e(t)
 \right\|_{L^2(\R^n)}\Tend \Lambda
{+\infty} 0\, ,\\
\text{where}\quad  &v^\e(t,x)=\frac{e^{i\frac{|x|^2}{2\e(t-1)}}}{(1-t)^{n/2}}
 a_0\(\frac{x}{1-t}\)  
\exp\(i\e^{\g-1}\phi\(
  \frac{\e^\g}{1-t}\virgp\frac{x}{1-t}\)\).
\end{align*}
\end{theo}
We now comment this result. In the linear case $f\equiv 0$, the above
result holds with $\g =1$ and $\phi\equiv 0$ (see
Section~\ref{sec:free}). We recall in Section~\ref{sec:critical} that in the
critical case ``nonlinear caustic, linear WKB'', the same asymptotic as in the
linear case holds for $1-t\gg \e$. The case
$k<n$ is super-critical as far as nonlinear effects near $t=1$ are
concerned. We emphasize two important
features in the above result: the analysis stops sooner than
$1-t\gg \e$, and nonlinear effects cause the presence of the
(nontrivial) phase $\phi$. For $1-t\gg \e^\g$, we have 
\begin{equation*}
 \e^{\g-1}\phi\(
  \frac{\e^\g}{1-t},\frac{x}{1-t}\) \sim  \sum_{j\ge 1}
  \frac{\e^{jk-1}}{(1-t)^{jn-1}} \phi_j\(\frac{x}{1-t}\)\, .
\end{equation*}
The above phase shift starts being relevant for $1-t \sim
\e^{\frac{k-1}{n-1}}$ (recall that $n>k>1$); this is the first
boundary layer where nonlinear effects appear at leading
order, measured by $\phi_1$. We will check that this phase shift is
relevant: $\phi_1$ is not 
zero (unless $u^\e\equiv 0$, see \eqref{eq:phi1} below). 
We then have a countable number of boundary layers in time, of
size
\begin{equation*}
  1-t \sim
\e^{\frac{jk-1}{jn-1}}\, ,
\end{equation*}
which reach the layer $1-t\sim \e^\g$ in the limit $j\to
+\infty$. At each new boundary layer, a new phase $\phi_j$ becomes
relevant at leading order. In general, none of the $\phi_j$'s is zero:
see e.g. \eqref{eq:phi2} for $\phi_2$. 
The result of a cascade of phases can be compared to  
the one discovered recently by C.~Cheverry \cite{CheverryBullSMF} in
the case of fluid dynamics, although the phenomenon seems to be
different. Yet, our result shares another property with
\cite{CheverryBullSMF}, which does not appear in the above
statement. Theorem~\ref{theo:main} shows perturbations of the phase
(the $\phi_j$'s), but not of the amplitude: the main profile is the
same as in the linear case, that is, a rescaling of $a_0$. However, to
compute the first $N$ phase shifts, $(\phi_j)_{1\le j\le N}$, one has
to compute $N-1$ corrector terms of the main profile $a_0$. This appears in
Proposition~\ref{prop:small}; see also
Equations~\eqref{eq:phi1}--\eqref{eq:phi2}. 

The assumption $k>1$ means that we start with a linear WKB
r\'egime. Indeed, for small positive time, $u^\e$ remains of order
$\O(1)$, and $f(\e^k|u^\e|^2)\sim \e^k|u^\e|^2 f'(0)$. The main term
is then the same as in \cite{Ca2} with $\si=1$ and $\alpha =k$. As
recalled in the above table, $\alpha >1$ corresponds to a propagation
which is linear at leading order.

Each phase shift oscillates at a rate
between $\O(1)$ (when it starts being relevant) and $\O(\e^{\g-1})$ (when
it reaches the layer of size $\e^\g$). Since $\g>0$, this means that
each phase shift is rapidly oscillating at the scale of the amplitude,
but oscillating strictly more slowly than the geometric phase
$\frac{|x|^2}{2\e(t-1)}$, for $1-t\gg \e^\g$. We will see in
Section~\ref{sec:beyond} that for $1-t =\O(\e^\g)$, all the terms in
$\phi$, plus the geometric phase, have the same order: all these
phases become comparable, see \eqref{eq:chg}. 

We will prove a more precise asymptotics than the $L^2$ estimate of
Theorem~\ref{theo:main}: see Proposition~\ref{prop:bkwreduit} and
\eqref{eq:pseudo}. We restricted our attention to the $L^2$ norm for
the sake of brevity.
 
\begin{figure}[htbp]
\begin{center}
\input{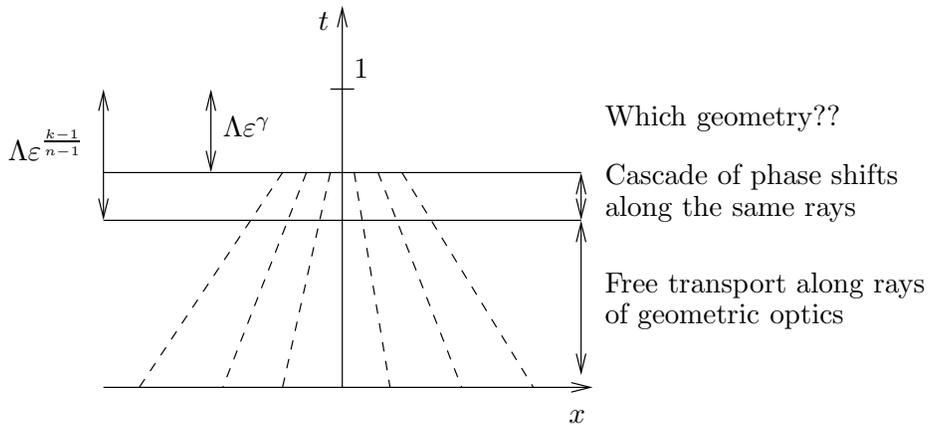}
\caption{Geometry of the propagation, with morally $\Lambda \to +\infty$.}
\label{fig:rays}
\end{center}
\end{figure}

Unfortunately, our analysis stops at the boundary layer of size
$\e^\g$: we can only go up to $1-t = \e^\g/T$ where $T$ appears in
Theorem~\ref{theo:main}. We will discuss this fact in
Section~\ref{sec:beyond}, and 
explain why we took care of never speaking of ``focal point'' in the
super-critical case, but only of caustic (as a matter of fact, even
the existence of a caustic is not clear, see
Section~\ref{sec:beyond}). For instance, the geometry 
of the propagation is not known for $1-t\lesssim \e^\g$, while the
analysis shows that it occurs on the rays of linear geometric optics
before this layer (see Figure~\ref{fig:rays}). 
On the other hand, we know that the order of magnitude of the
amplitude changes to reach the boundary layer of size $\e^\g$. Recall that
Theorem~\ref{theo:main} describes the asymptotic behavior of $u^\e$
for $1-t\ge \Lambda \e^\g$, in the limit $\Lambda \to +\infty$. In this
region, leading order nonlinear effects are visible only in the phase. As
mentioned above, our analysis is valid for $1-t \ge  \e^\g/T$.
Between the initial time and this region, the amplitude
of $u^\e$ varies like $(1-t)^{-n/2}$, and changes from $\O(1)$ to
$\O(\e^{-k/2})$ (see also Remark~\ref{rema:bkwameliore} for the
asymptotic behavior of $u^\e$ for $t\le 1-\e^\g/T$).

The rest of the paper is organized as follows. In
Section~\ref{sec:rappel}, we recall the proof of the analog of
Theorem~\ref{theo:main} in the linear and critical nonlinear cases. In
Section~\ref{sec:formal}, we present a formal computation that suggests
a result like Theorem~\ref{theo:main}. Based on the
result by E.~Grenier
\cite{Grenier98} and a ``semi-classical conformal transform'',
we give the proof of Theorem~\ref{theo:main} in
Section~\ref{sec:justif}. In Section~\ref{sec:stab}, we compare the
rigorous approach with 
the formal result of Section~\ref{sec:formal}. The
discussion about some possible phenomena for $t\ge 1- \e^\g/T$
appears in the final Section~\ref{sec:beyond}. 

\bigskip

\noindent 
{\bf Acknowledgments. } The author is grateful to Christophe Cheverry 
for stimulating discussions about this work. 
\section{Free and critical cases}
\label{sec:rappel}

\subsection{The linear equation}
\label{sec:free}
Consider the linear equation:
\begin{equation}
  \label{eq:lin}
  i\e \d_t u_{\rm lin}^\e + \frac{\e^2}{2} \Delta u_{\rm lin}^\e
  =0\, ,\quad (t,x)\in 
  \R\times \R^n\ ; \quad
  u^\e_{{\rm lin}\mid t=0}= a_0(x) e^{-i\frac{|x|^2}{2\e}} .  
\end{equation}
As $\e \to 0$, the rays of geometric optics (classical trajectories)
are lines $\frac{x}{1-t}=\text{const.}$, and meet at the origin at time
$t=1$. Indeed, the bicharacteristic curves are
defined by the Hamilton flow associated to $p(t,x,\tau,\xi) = \tau
+\frac{|\xi|^2}{2}$: 
\begin{equation*}
  \dot t=1\quad ;\quad \dot x=\xi\quad ;\quad \dot \tau =\dot \xi
  =0\quad ; \quad x(0)=x_0 \quad ;\quad \xi(0)= \nabla\phi(0,x(0))= -x_0\, . 
\end{equation*}
Of course, $u_{\rm 
  lin}^\e$ can be expressed in terms of an oscillatory
integral:
\begin{equation}\label{eq:sollin}
  u^\e(t,x) = \frac{1}{(2i\pi \e t)^{n/2}}\int_{\R^n}
  e^{i\frac{|x-y|^2}{2\e t} -i\frac{|y|^2}{2\e}}a_0(y)dy\, .
\end{equation}
Applying stationary phase formula yields the same result as
using WKB methods below, up to the same boundary layer. Seek
\begin{equation*}
  u_{\rm lin}^\e (t,x)\Eq \e 0 v_{\rm lin}^\e (t,x) =
  v^0(t,x)e^{i\frac{\phi(t,x)}{\e}}\, .  
\end{equation*}
Plugging this into \eqref{eq:lin} and canceling the $\O(\e^0)$ and
$\O(\e^1)$ terms, we find:
\begin{align*}
  \d_t \phi +\frac{1}{2}|\nabla_x \phi|^2 &= 0\, ,\quad
  \phi(0,x)=-\frac{|x|^2}{2}\ ; \\
\d_t v^0 + \nabla_x \phi\cdot \nabla_x v^0 +\frac{1}{2}v^0 \Delta \phi
&=0 \, , \quad v^0(0,x)=a_0(x)\, .
\end{align*}
For $t<1$, one has explicitly:
\begin{equation*}
  \phi(t,x)= \frac{|x|^2}{2(t-1)}\quad ; \quad
  v^0(t,x)=\frac{1}{(1-t)^{n/2}} a_0\(\frac{x}{1-t}\)\, .
\end{equation*}
Moreover, $v_{\rm
  lin}^\e$ solves: 
\begin{equation*}
  i\e \d_t v_{\rm lin}^\e + \frac{\e^2}{2} \Delta v_{\rm lin}^\e =
  \frac{\e^2}{2(1-t)^2}\frac{e^{i\frac{|x|^2}{2(t-1)}} }{(1-t)^{n/2}}
  \Delta a_0 \(\frac{x}{1-t}\)
  \   ; \ 
  v^\e_{{\rm lin}\mid t=0}= a_0(x) e^{-i\frac{|x|^2}{2\e}} .  
\end{equation*}
Let $r^\e(t,x)$ denote the source term:
$\dis   \|r^\e(t)\|_{L^2}\lesssim \frac{\e^2}{(1-t)^2}$. 
Standard energy estimates for Schr\"odinger equation yield:
\begin{equation*}
  \e \frac{d}{dt}\| u^\e(t) - v^\e(t)\|_{L^2} \lesssim
  \|r^\e(t)\|_{L^2}\lesssim \frac{\e^2}{(1-t)^2}\, ,
\end{equation*}
and 
\begin{equation*}
  \sup_{0\le s\le t}\| u^\e(s) - v^\e(s)\|_{L^2} \lesssim
  \frac{1}{\e}\int_0^t\|r^\e(s)\|_{L^2}ds \lesssim \frac{\e}{1-t}\, \cdot
\end{equation*}
Thus, WKB approximation is interesting up to a boundary
layer in time of size $\e$ before the focus. Moreover, $v_{\rm
  lin}^\e$ is exactly the approximate solution of 
Theorem~\ref{theo:main} with $\phi\equiv 0$. 
Past this boundary layer, \eqref{eq:sollin} shows that for
$|1-t|=\O(\e)$, 
\begin{equation*}
  u^\e(t,x) \sim \frac{e^{i\frac{|x|^2}{2\e}}}{(2i\pi \e)^{n/2}} \int
  e^{-i\frac{x\cdot y}{\e}}a_0(y)dy =
  \frac{e^{i\frac{|x|^2}{2\e}}}{\e^{n/2}} \F a_0\(\frac{x}{\e}\)\sim
  \frac{1}{\e^{n/2}} \F a_0\(\frac{x}{\e}\) ,
\end{equation*}
where $\F$ denotes the Fourier transform.
For $t-1\gg \e$, stationary phase formula yields 
the same asymptotic description as above, up to the Maslov
index (see \cite{Du,Ca2}). In particular, we see that
\eqref{eq:estlinsharp} holds for $u^\e_{\rm lin}$, and is sharp. 
\subsection{The critical nonlinear case}
\label{sec:critical}
We recall the main result of \cite{Ca2}. 
Consider \eqref{eq:princ} in the case $k=n$, which is critical
concerning the role of the nonlinearity near the focal point. 
Introduce the scaling
\begin{equation*}
  u^\e(t,x)=\frac{1}{\e^{n/2}}\psi^\e\( \frac{t-1}{\e}\virgp
  \frac{x}{\e}\) .
\end{equation*}
Then the function $\psi^\e$ solves 
\begin{equation}
  i \d_t \psi^\e + \frac{1}{2} \Delta \psi^\e
  =f\(|\psi^\e|^{2}\)\psi^\e \quad ;\quad \psi^\e_{\mid t =
  \frac{-1}{\e}} = \e^{n/2}a_0(\e x)e^{-i\e\frac{|x|^2}{2}}\, . 
\end{equation}
A way to understand criticality is that $\e$ has disappeared from the
equation satisfied by $\psi^\e$. 
Using global well-posedness results for nonlinear Schr\"odinger
equations (under assumptions on the nonlinearity which are different
from Assumptions~\ref{hyp:tout}, see e.g. \cite{CazCourant}), one has
\begin{equation*}
  \left\| \psi^\e -\psi\right\|_{L^\infty(\R;H^1(\R^n))}\to 0 \text{
  as }\e \to 0\, ,
\end{equation*}
where $\psi$ is the (global) solution of the Cauchy problem
\begin{equation}\label{eq:scatt}
  i \d_t \psi + \frac{1}{2} \Delta \psi =f\(|\psi|^{2}\)\psi\quad ;
  \quad e^{-i\frac{t}{2}\Delta}\psi(t,x)\big|_{t=-\infty} =
  \F^{-1}(a_0)(x)\, . 
\end{equation}
Scattering theory shows two interesting features: for large
$|t|$, $\psi(t,x)$ behaves like a solution of the linear Schr\"odinger
equation. This implies that for $1-t\gg \e$, the solution $u^\e$ can be
approximated by $u_{\rm lin}^\e$ (or $v_{\rm lin}^\e$): no nonlinear
effect is relevant before the same boundary layer as before. The
second point is that for $|1-t|\lesssim \e$, nonlinear effects occur
at leading order, and are measured (in average) by the nonlinear
scattering operator associated to \eqref{eq:scatt}. 

\section{A formal computation}
\label{sec:formal}

From now on, we assume that $k<n$. To simplify notations, and since
the Assumptions~\ref{hyp:tout} will be needed for rigorous proofs only,
consider the case of an homogeneous nonlinearity:
$  f(y) =y^\si$, and denote $\alpha = k\si$. Then \eqref{eq:princ}
becomes
\begin{equation}
  \label{eq:homo}
  i\e \d_t u^\e + \frac{\e^2}{2} \Delta u^\e = \e^\alpha
  |u^\e|^{2\si} u^\e \quad ; \quad 
  u^\e_{\mid t=0}= a_0(x) e^{-i\frac{|x|^2}{2\e}} .  
\end{equation}
The caustic is supercritical: $n\si>\alpha$. We also assume $\alpha>1$
(linear WKB).  Because this section is only formal, we shall be very
brief about the computations, and only give the main steps. 
\subsection{A first boundary layer}
Two approaches (at least) lead to the same result: Lagrangian integral
with stationary phase formula (like in \cite{CaCMP} where the critical
one-dimensional cubic case is considered), and generalized WKB
methods. We shall retain the second one, which we use in the next
subsection.  Seek
\begin{equation*}
  u^\e (t,x)\Eq \e 0 v^\e_1 (t,x) =
  u^0(t,x)e^{i\frac{\phi(t,x)}{\e}}\, ,
\end{equation*}
and change the usual hierarchy to force the contribution of the
nonlinear term to appear in the transport equation:
\begin{align*}
  \d_t \phi +\frac{1}{2}|\nabla_x \phi|^2 &= 0\, ,\quad
  \phi(0,x)=-\frac{|x|^2}{2}\ ; \\
\d_t u^0 + \nabla_x \phi\cdot \nabla_x u^0 +\frac{1}{2}u^0 \Delta \phi
&=-i \e^{\alpha -1}|u^0|^{2\si}u^0  \, , \quad u^0(0,x)=a_0(x)\, .
\end{align*}
The eikonal equation is the same as in Section~\ref{sec:free}, as
well as its solution. 
The transport equation is an ordinary differential equation along the
rays of geometric optics $\frac{x}{1-t}=\text{const.}$, of the form
\begin{equation*}
  \dot y = -i \e^{\alpha -1}|y|^{2\si}y\, .
\end{equation*}
The modulus of  $u^0$ is constant along rays, and 
\begin{equation*}
 u^0(t,x) = \frac{1}{(1-t)^{n/2}} a_0\(\frac{x}{1-t}\)
 \exp\(-i\e^{\alpha -1}\left|
 a_0\(\frac{x}{1-t}\)\right|^{2\si}\int_0^t\frac{ds}{(1-s)^{n\si}}\).  
\end{equation*}
Note that the notation is no longer relevant, since $u^0$ now depends
on $\e$. 
We have a new boundary layer in time, of size $\e^\beta$ before the
focus, where 
\begin{equation*}
  \beta = \frac{\alpha -1}{n\si -1}\, .
\end{equation*}
For $1-t \sim \e^\beta$, the above phase shift measures relevant
nonlinear effects. We have:
\begin{equation*}
  i\e \d_t v^\e_1 + \frac{\e^2}{2} \Delta v^\e_1 = \e^\alpha
  |v^\e_1|^{2\si}v^\e_1 + r_1^\e \quad ; \quad
  v^\e_{1\mid t=0}= a_0(x) e^{-i\frac{|x|^2}{2\e}}\, ,  
\end{equation*}
with
\begin{equation*}
  \frac{1}{\e}\int_0^t\|r_1 ^\e(s)\|_{L^2}ds \lesssim \frac{\e}{1-t} +
  \frac{\e^{2\alpha -1}}{(1-t)^{2n\si -1}}\, .
\end{equation*}
Following the energy estimates of Section~\ref{sec:free}, this
quantity might be the one that dictates the size of the error $u^\e
-v^\e_1$ (see Section~\ref{sec:stab} for a discussion on that issue).  
The second term is ``new'' (the first term is the same as in
Section~\ref{sec:free}), and suggests the existence of a second
boundary layer, of size
$  \e^{\frac{2\alpha -1}{2n\si -1}}$. 
\subsection{Infinitely many boundary layers: cascade of phase shifts}
Seek an approximate solution of the form:
\begin{equation*}
  v^\e (t,x) = \frac{1}{(1-t)^{n/2}} a_0\(\frac{x}{1-t}\) e^{i
    \phi^\e(t,x)},\quad \phi^\e(t,x) = \frac{|x|^2}{2\e(t-1)} + g^\e
  (t,x)\, .
\end{equation*}
We find
\begin{equation*}
\begin{aligned}
  i\e \d_t v^\e + \frac{\e^2}{2}& \Delta v^\e =  \( i\frac{\e^2}{2} 
\Delta g^\e
  -\e \d_t g^\e -\frac{\e^2}{2}|\nabla_x g^\e|^2 + \frac{\e}{1-t}
  x\cdot \nabla_x g^\e \) v^\e\\
 + i\frac{\e^2}{(1-t)^{\frac{n}{2}+1}}&\nabla_x g^\e\cdot \nabla a_0
\(\frac{x}{1-t}\) e^{i  \phi^\e} + \frac{1}{2}\( \frac{\e}{1-t}\)^2
\frac{e^{i  \phi^\e}}{(1-t)^{n/2}}\Delta a_0 \(\frac{x}{1-t}\)\, .
\end{aligned} 
\end{equation*}
As suggested by the previous paragraph, write 
\begin{equation}
  \label{eq:g}
  g^\e (t,x) = \frac{1}{\e} \int_{0}^t h\(
  \frac{\e^\alpha}{(1-s)^{n\si}}\virgp \frac{x}{1-t}\) ds\, ,\quad
  \text{with } h(z,\xi) \sim \sum_{j\ge 1}z^j g_j(\xi)\, .
\end{equation}
In the equation solved by $v^\e$, the last term is the ``same'' as in
the linear case: it becomes relevant only in a boundary layer of size
$\e$. Since our approach will lead us to the boundary layer of size
$\e^\g$ (recall that $\g =k/n = \alpha/n\si<1$), we ignore that
term.  

The remaining terms with a factor $i$ are of order, in $L^2$, 
\begin{equation*}
\e^2\|\Delta g^\e(t)\|_{L^\infty} + \frac{\e^2}{1-t}\|\nabla_x
g^\e(t)\|_{L^\infty} \lesssim \frac{\e}{(1-t)^2}\int_{t_0}^t
\frac{\e^\alpha}{(1-s)^{n\si}} ds \lesssim
\frac{\e^{\alpha+1}}{(1-t)^{n\si+1}}\, , 
\end{equation*}
and their contribution is also left out in this computation. 

Now we require that $v^\e$ be an approximate solution to \eqref{eq:homo}:
\begin{equation}
  \label{eq:evolg}
  \( \d_t -\frac{x}{1-t}\cdot \nabla_x\) g^\e + \frac{\e}{2}|\nabla_x
  g^\e|^2 =  -\frac{\e^{\alpha -1}}{(1-t)^{n\si}}\left|
    a_0\(\frac{x}{1-t}\)\right|^{2\si}\, . 
\end{equation}
Using \eqref{eq:g}, we get:
\begin{equation}
  \label{eq:solg}
  \begin{aligned}
    & g_1(\xi)= -|a_0(\xi)|^{2\si}\, ,\\
\text{for }j\ge 2, \quad & g_j (\xi) = -\frac{1}{2} \sum_{p+q=j}
\frac{1}{(pn\si -1)(qn\si -1)} \nabla g_p\cdot \nabla g_q\, ,
  \end{aligned}
\end{equation}
with the convention $g_0\equiv 0$. This algorithm produces smooth
solutions provided that $|a_0(\xi)|^{2\si}$ is smooth ($\si\in \N^*$
or $a_0$ Gaussian for instance). We neglected the terms corresponding
to $s=0$ in the integration \eqref{eq:g}: this does not increase the
error, since  
$n\si >\alpha > 1$. Defining 
\begin{align*}
   \widetilde g_N^\e (t,x) &= \frac{1}{\e} \sum_{j=1}^N \int_0^t \(
  \frac{\e^\g}{1-s} \)^{n\si j} ds\times g_j\(\frac{x}{1-t}\),\\
\widetilde v^\e_N (t,x)& = \frac{1}{(1-t)^{n/2}} a_0\(\frac{x}{1-t}\)
  e^{i\frac{|x|^2}{2\e(t-1)} + \widetilde g_N^\e 
  (t,x)}\, ,
\end{align*}
the approximate solution $\widetilde v_N^\e$ solves 
\begin{equation*}
  i\e \d_t \widetilde v^\e_N + \frac{\e^2}{2} \Delta \widetilde v^\e_N =
  \e^\alpha 
  |\widetilde v^\e_N|^{2\si}\widetilde v^\e_N + \widetilde r_N^\e
  \quad ;\quad 
  \widetilde v^\e_{N\mid t=0}=  a_0(x)e^{-i\frac{|x|^2}{2\e}}\, ,
\end{equation*}
with, for $1-t\ge \e^\g$:
\begin{equation*}
 \frac{1}{\e}\int_0^t \| \widetilde  r^\e_N(s)\|_{L^2}ds \lesssim
 \frac{\e^{(N+1)\alpha -1}}{(1-t)^{(N+1)n\si-1}}+ 
 \frac{\e^{\alpha}}{(1-t)^{n\si}} \, .
\end{equation*} 
To compare with Theorem~\ref{theo:main},
remove the terms corresponding to $s=0$ in the integration (recall
that $n\si>1$), and define:
\begin{align*}
   g_N^\e (t,x) &= \frac{1}{\e} \sum_{j=1}^N \int_{-\infty}^t \(
  \frac{\e^\g}{1-s} \)^{n\si j} ds\times g_j\(\frac{x}{1-t}\),\\
v^\e_N (t,x) &= \frac{1}{(1-t)^{n/2}} a_0\(\frac{x}{1-t}\)
  e^{i\frac{|x|^2}{2\e(t-1)} + g_N^\e 
  (t,x)}\, .
\end{align*}
By definition, we have $\|v_N^\e (t)-\widetilde v_N^\e
(t)\|_{L^2} = \O(\e^{\alpha -1})$ for $1-t\ge \e^{\g}$. 
One can check that $v_N^\e$ solves 
\begin{equation*}
  i\e \d_t v^\e_N + \frac{\e^2}{2} \Delta v^\e_N =
  \e^\alpha 
  |v^\e_N|^{2\si}v^\e_N + r_N^\e \, , 
\end{equation*}
with
\begin{equation}
  \label{eq:r_N}
\begin{aligned}
  r_N^\e (t,x)=&  \(\frac{\e^\alpha}{(1-t)^{n\si}}\)^{N+1}
  v_N^\e (t,x) g_{N+1}\(\frac{x}{1-t}\)\\
  +\frac{1}{2}&\(\frac{\e}{1-t}\)^2 e^{ig_N^\e 
  (t,x)} \frac{e^{i\frac{|x|^2}{2\e(t-1)}}}{(1-t)^{n/2}}\Delta a_0
  \(\frac{x}{1-t}\)   \\
+ i\frac{\e}{1-t} &\sum_{j=1}^N \frac{1}{jn\si-1}
  \(\frac{\e^\alpha}{(1-t)^{n\si}}\)^j \bigg( v_N^\e(t,x)\Delta g_j
  \(\frac{x}{1-t}\) \\  
&\ \ \ \ \ \ \ \ \ \ \ \  +e^{ig_N^\e 
  (t,x)} \frac{e^{i\frac{|x|^2}{2\e(t-1)}}}{(1-t)^{n/2}} \nabla g_j
  \cdot \nabla a_0\(\frac{x}{1-t}\)\bigg). 
\end{aligned}
\end{equation}
We have the following result:
\begin{prop}[Formal approximation to \eqref{eq:homo}]\label{prop:formal}
  Let $n\si>\alpha>1$, $a_0\in \S(\R^n)$, and fix $N\in \N^*$. Denote 
\begin{equation*}
  g_N^\e (t,x) =\sum_{j=1}^N
  \frac{\e^{j\alpha-1}}{(1-t)^{jn\si -1}} \frac{1}{jn\si-1}
  g_j\(\frac{x}{1-t}\) , 
\end{equation*}
and let $v_N^\e$ be the associated approximate solution. The
  function  $v_N^\e$ solves
\begin{equation*}
  i\e \d_t v^\e_N + \frac{\e^2}{2} \Delta v^\e_N =
  \e^\alpha 
  |v^\e_N|^{2\si}v^\e_N + r_N^\e \  ; \  v^\e_{N\mid t=0} =
  a_0(x)e^{-i\frac{|x|^2}{2\e}} + \O\(\e^{\alpha -1}\) \text{ in }L^2.   
\end{equation*}
For $1-t \ge \e^{\g}=\e^{\frac{\alpha}{n\si}}$, the source term
satisfies: 
\begin{equation*}
 \frac{1}{\e}\int_0^t \|  r^\e_N(s)\|_{L^2}ds \lesssim
 \frac{\e^{(N+1)\alpha -1}}{(1-t)^{(N+1)n\si-1}}+
 \frac{\e^{\alpha}}{(1-t)^{n\si}} \, .
\end{equation*} 
\end{prop}
For
$1\le j\le N$, the $j^{\rm th}$ term of the series defining $g_N^\e$
becomes relevant in a 
boundary layer of size $\e^{\frac{j\alpha -1}{nj\si -1}}$: in the
limit $N\to +\infty$, a countable
family of boundary 
layers appear, between $\e^\beta$ and $\e^\g$. In the case $\si =1$,
which is the only homogeneous nonlinearity consistent with
Assumptions~\ref{hyp:tout}, 
we have $\alpha =k$ and we find the boundary layers announced in the
introduction. 

Letting $N\to +\infty$ (using Borel lemma, see e.g. \cite{RauchUtah}),
we find: 
\begin{equation*}
 \frac{1}{\e}\int_0^t \|  r^\e(s)\|_{L^2}ds \lesssim
 \frac{\e^{\alpha}}{(1-t)^{n\si}} \, ,
\end{equation*}
which is small for $1-t\gg \e^\gamma$. 
\begin{rema}
  In the critical case $\alpha =n\si >1$, we have $\beta = \g =1$:
  the above boundary layers ``collapse'' one on another. There are no
  such phase shifts as above. 
\end{rema}
We point out that the sole estimate of the source term proves
nothing. In a stability argument, the nonlinearity $|u^\e|^{2\si}u^\e-
|v^\e_N|^{2\si}v^\e_N$ is usually treated by a Gronwall type
argument. If the nonlinearity is  ``too strong'', then the above
estimate, which is completely relevant in the linear case, does not
necessarily account for the size of the error. Since we are in a
super-critical 
case, it is not surprising that Proposition~\ref{prop:formal} is only
a formal result. This remark can be compared to the approach in
\cite{PGX93}. To justify a WKB expansion for the nonlinear equation
\begin{equation*}
  i\e \d_t u^\e +\frac{\e^2}{2}\Delta u^\e = f\(|u^\e|^2\)u^\e\, ,
\end{equation*}
constructing an approximate solution that solves
\begin{equation*}
  i\e \d_t u^\e_{\rm app} +\frac{\e^2}{2}\Delta u^\e_{\rm app} =
  f\(|u^\e_{\rm app}|^2\)u^\e_{\rm app}+\O\(\e^\infty\)\, , 
\end{equation*}
is not sufficient. Indeed, the computations in \cite{PGX93} show that
energy estimates and Gronwall lemma do not yield better than
\begin{equation*}
  \left\|u^\e(t)-u^\e_{\rm app}(t)\right\|_{L^2}\le
  e^{Ct/\e}\O\(\e^\infty\). 
\end{equation*}
This is the reason why in \cite{PGX93}, WKB expansions are justified
for analytic data. This assumption yields a source term for $u^\e_{\rm
  app}$ which is $\O(e^{-\delta/\e})$, counterbalancing the
exponential growth of Gronwall lemma to lead to a good approximation
on $[0,T]$ for $T>0$ independent of $\e$. 

\section{Rigorous results}
\label{sec:justif}

We now prove Theorem~\ref{theo:main}. We will see that the approximate
solution we find diverges from the one constructed above, a fact which
we discuss in Section~\ref{sec:stab}.

\subsection{Semi-classical conformal transform}
\label{sec:pseudo}
Introduce the new unknown function $\psi^\e$ given by:
\begin{equation}
  \label{eq:pseudo}
  u^\e(t,x) =
  \frac{1}{(1-t)^{n/2}}\psi^\e\(\frac{\e^\g}{1-t}\virgp
  \frac{x}{1-t} \) e^{i\frac{|x|^2}{2\e(t-1)}}\, .
\end{equation}
Recalling that $\g=\frac{k}{n}<1$, denote
\begin{equation}
  \label{eq:h}
  \h = \e^{1-\g}\to 0\, .
\end{equation}
Changing the notation $\psi^\e(\tau,\xi)$ into $\psi^\h(t,x)$, we check that
\eqref{eq:princ} becomes:
\begin{equation}
  \label{eq:psi}
  i\h \d_t \psi^\h +\frac{\h^2}{2}\Delta \psi^\h = t^{-2}f\( t^n |\psi^\h|^2\)
  \psi^\h\quad ; \quad \psi^\h\big|_{ t = \h^{\frac{\g}{1-\g}}} = a_0(x)\, .
\end{equation}
The singular term $t^{-2}$ in factor of the nonlinearity is actually
harmless: as $t$ goes to zero, $t^{-2} f\( t^n
|\psi^\h|^2\)\thickapprox t^{n-2} 
|\psi^\h|^2 f'(0)$, and is bounded since $n\ge 2$.

The proof of  Theorem~\ref{theo:main} is now reduced to the asymptotic
expansion for $\psi^\h$ as $\h\to 0$ for $t\in [\h^{\frac{\g}{1-\g}},
\frac{1}{\Lambda}]$. Denote $t_0^\h= \h^{\frac{\g}{1-\g}}$. We shall
prove the following: 
\begin{prop}
\label{prop:bkwreduit}
Let Assumptions~\ref{hyp:tout} be satisfied. Assume $n>k>1$, and let
$s\in \N$. Then there exists $T>0$ independent of $\h\in ]0,1]$ such
that for $t\in [t_0^\h,T]$, the function $\psi^\h$ can be written as 
$\psi^\h(t,x)=a^\h(t,x)e^{i\phi^\h(t,x)/\h}$, with
\begin{equation*}
  \left\|a^\h-a\right\|_{L^\infty([t_0^\h,T];H^s)}
  +\left\|\phi^\h-\phi\right\|_{L^\infty([t_0^\h,T];H^s)}  \to
  0 \quad \text{as }\h \to 0\, ,
\end{equation*}
where $(a,\phi)$ solves
\begin{equation}
  \label{eq:bkwlimite}
  \begin{aligned}
    \d_t \phi +\frac{1}{2}|\nabla_x \phi|^2 +t^{-2}f\(t^n |a|^2\)
    &=0\quad ; \quad \phi_{\mid t=0}=0\, ,\\
\d_t a +\nabla_x \phi \cdot \nabla_x a +\frac{1}{2}a\Delta \phi
    &= 0\quad ; \quad a_{\mid t=0}=a_0\, .
  \end{aligned}
\end{equation}
Moreover,
\begin{equation*}
 \limsup_{\h \to 0}\sup_{t_0^\h\le t\le \tau}\left\| \psi^\h(t,x)
 -a_0\(x\) 
e^{i\phi(t,x)/\h }\right\|_{H^s(\R^n)}\to 0 \quad \text{as }\tau \to 0\, .
\end{equation*}
\end{prop}
The second point of Theorem~\ref{theo:main} follows from the above
proposition, since  
the transform \eqref{eq:pseudo} is $L^2$ unitary (see
Proposition~\ref{prop:small} below for the asymptotic expansion of
$\phi$). We could also
include not only derivatives in the above estimates, but also
momenta. As announced in the introduction, we chose to leave out this
refinement. Note
that except for two 
aspects, Proposition~\ref{prop:bkwreduit} is nothing but rewriting
Theorems~1.1 and 1.3 of \cite{Grenier98}. In our case, time is present
in the nonlinearity, and data for $\psi^\h$ are prescribed at time
$t_0^\h$ (with $t_0^\h\to 0$ as $\h \to 0$) instead of time zero.  
\subsection{Construction of solutions to \eqref{eq:psi}}
\label{sec:grenier}

We recall the ideas introduced by E.~Grenier
\cite{Grenier98}, and show how to handle the presence of time in the
nonlinearity. The main idea in \cite{Grenier98} is to write the
solution of a semi-classical nonlinear Schr\"odinger equation as a WKB
solution, where not only the amplitude may depend on the small
parameter, but also the phase. This changes the usual WKB hierarchy,
and overcomes the difficulties pointed out in \cite{PGX93}. Seek the
solution of \eqref{eq:psi} of the form
\begin{equation*}
  \psi^\h(t,x) = a^\h (t,x)e^{i\phi^\h(t,x)/\h},
\end{equation*}
with
\begin{equation}\label{eq:systexact0}
  \begin{aligned}
    \d_t \phi^\h +\frac{1}{2}\left|\nabla \phi^\h\right|^2 + t^{-2}f\(t^n
    |a^\h|^2\)= 0\quad &; \quad
    \phi^\h\big|_{t=t_0^\h}=0\, ,\\
\d_t a^\h +\nabla \phi^\h \cdot \nabla a^\h +\frac{1}{2}a^\h
\Delta \phi^\h  = i\frac{\h}{2}\Delta a^h\quad & ;\quad
a^\h\big|_{t=t_0^\h}= a_0\, . 
  \end{aligned}
\end{equation}
Introducing the ``velocity'' ${\tt v}^\h = \nabla \phi^\h$,
\eqref{eq:systexact0} yields
\begin{equation}\label{eq:systexact}
  \begin{aligned}
    \d_t {\tt v}^\h +{\tt v}^\h \cdot \nabla {\tt v}^\h + 2 t^{n-2}f'\(t^n
    |a^\h|^2\) \operatorname{Re}\(\overline{a^\h}\nabla a^\h\)=
    0\quad &; \quad 
    {\tt v}^\h\big|_{t=t_0^\h}=0\, ,\\
\d_t a^\h +{\tt v}^\h\cdot \nabla a^\h +\frac{1}{2}a^\h
\operatorname{div}{\tt v}^\h  = i\frac{\h}{2}\Delta a^h\quad & ;\quad
a^\h\big|_{t=t_0^\h}= a_0\, . 
  \end{aligned}
\end{equation}
To force the initial time to be zero, introduce
\begin{equation*}
  \widetilde {\tt v}^\h(t,x)={\tt v}^\h\(t+t_0^\h,x\)\quad ; \quad
  \widetilde a^\h(t,x)=a^\h\(t+t_0^\h,x\) . 
\end{equation*}
Then \eqref{eq:systexact} becomes
\begin{equation*}
  \begin{aligned}
    \d_t \widetilde {\tt v}^\h +\widetilde {\tt v}^\h \cdot \nabla
    \widetilde {\tt v}^\h + 2 \(t+t_0^\h\)^{n-2}f'\(\(t+t_0^\h\)^n 
    |\widetilde a^\h|^2\) \operatorname{Re}
\(\overline{\widetilde a^\h}\nabla \widetilde a^\h\)=
    0\  &; \ 
    \widetilde {\tt v}^\h\big|_{t=0}=0 ,\\
\d_t \widetilde a^\h +\widetilde {\tt v}^\h\cdot \nabla \widetilde
    a^\h +\frac{1}{2} \widetilde a^\h
\operatorname{div}\widetilde {\tt v}^\h  = i\frac{\h}{2}\Delta
    \widetilde a^h\  & ;\  
\widetilde a^\h\big|_{t=0}= a_0 . 
  \end{aligned}
\end{equation*}
Notice that if $n=2$ and $f'=\text{const.}$ (2D cubic equation, which
is conformally invariant), the above system is exactly the same as in
\cite{Grenier98}. 

Separate real and
imaginary parts of $\widetilde a^\h$, $\widetilde a^\h = \widetilde
{\tt a}^\h + i \widetilde {\tt b}^\h$. Then we have 
\begin{equation}
  \label{eq:systhyp}
  \d_t \bu^\h +\sum_{j=1}^n A_j(\bu^\h)\d_j \bu^\h
  = \frac{\h}{2} L 
  \bu^\h\, , 
\end{equation}
\begin{equation*}
  \text{with}\quad \bu^\h = \left(
    \begin{array}[l]{c}
      \widetilde {\tt a}^\h \\
      \widetilde {\tt b}^\h \\
      \widetilde {\tt v}^\h_1 \\
      \vdots \\
      \widetilde {\tt v}^\h_n
    \end{array}
\right)\quad , \quad L = \left(
  \begin{array}[l]{ccccc}
   0  &-\Delta &0& \dots & 0   \\
   \Delta  & 0 &0& \dots & 0  \\
   0& 0 &&0_{n\times n}& \\
   \end{array}
\right),
\end{equation*}
\begin{equation*}
  \text{and}\quad A(\bu,\xi)=\sum_{j=1}^n A_j(\bu)\xi_j
= \left(
    \begin{array}[l]{ccc}
     \widetilde {\tt v}\cdot \xi & 0& \frac{\widetilde {\tt
     a}}{2}\,^{t}\xi \\ 
     0 & \widetilde {\tt v}\cdot \xi & \frac{\widetilde {\tt
     b}}{2}\,^{t}\xi \\ 
     2\(t+t_0^\h\)^{n-2}f' \widetilde {\tt a} \, \xi
     &2\(t+t_0^\h\)^{n-2}f' \widetilde {\tt b}\, \xi &  \widetilde
     {\tt v}\cdot \xi I_n 
    \end{array}
\right),
\end{equation*}
where $f'$ stands for $f'\( \(t+t_0^\h\)^n(|\widetilde {\tt a}|^2 +
  |\widetilde {\tt 
  b}|^2)\)$. The matrix $A(\bu,\xi)$ can be symmetrized by 
\begin{equation*}
  S=\left(
    \begin{array}[l]{cc}
     I_2 & 0\\
     0& \frac{1}{4\(t+t_0^\h\)^{n-2}f'}I_n
    \end{array}
\right),
\end{equation*}
which is symmetric and positive since $f'>0$. We now reproduce the
ideas of \cite{Grenier98}, inspired by hyperbolic theory, see
e.g. \cite{Majda}.  For an integer $s>2+n/2$, we bound $(S
\d_{x}^{\alpha } \bu^\h 
, \d_{x}^{\alpha }\bu^\h )$ 
where $\alpha $ is a  multi index of length $\le s$, and
$(\cdot ,\cdot)$ is the usual  $L^{2}$ scalar product. We have
\begin{equation*}
\frac{d}{dt}\(S \partial _{x}^{\alpha } \bu^\h ,
  \partial _{x}^{\alpha } \bu^\h\) 
= \(\partial _{t} S  \partial _{x}^{\alpha } \bu^\h , \partial
_{x}^{\alpha } \bu^\h\) 
  + 2 \( S \partial _{t}  \partial _{x}^{\alpha } \bu^\h , \partial
  _{x}^{\alpha } \bu^\h\) 
 \end{equation*}
since $S$ is symmetric. For the first term, we must consider the lower
$n\times n$ block in $S$. Differentiating $(t+t_0^\h)^{2-n}$ yields a
non-positive term, and we get
\begin{equation*}
\(\partial _{t} S  \partial _{x}^{\alpha } \bu^\h ,
  \partial _{x}^{\alpha } \bu^\h\) 
\le  \left\|\frac{1}{f'}\d_t\(f'\( \(t+t_0^\h\)^n(|\widetilde {\tt
  a}^\h|^2 + |\widetilde {\tt 
  b}^\h|^2)\)\)\right\|_{L^\infty}\(S  \partial _{x}^{\alpha } \bu^\h ,
  \partial _{x}^{\alpha } \bu^\h\)\, .
\end{equation*}
So long as $\|\bu^\h\|_{L^\infty}\le 2\|a_0\|_{L^\infty}$, we have,
for $t\le 2$ (to fix the ideas),
\begin{equation*}
  f'\( \(t+t_0^\h\)^n(|\widetilde {\tt a}^\h|^2 + |\widetilde {\tt
  b}^\h|^2)\) \ge \inf\left\{ f'(y)\ ;\ 0\le y\le
  2^{n+2}\|a_0\|_{L^\infty}^2\right\} =\delta_n>0\, ,
\end{equation*}
where $\delta_n$ is now fixed, since $f'$ is continuous with
$f'>0$. We infer, for $t\le 2$,  
\begin{equation*}
 \left\|\frac{1}{f'}\d_t\(f'\( \(t+t_0^\h\)^n(|\widetilde {\tt
  a}^\h|^2 + |\widetilde {\tt 
  b}^\h|^2)\)\)\right\|_{L^\infty}\le C\(\|\bu^\h
  \|_{L^\infty}\)\|\d_t \bu^\h 
  \|_{L^\infty} \lesssim \|\bu^\h\|_{H^s}\, ,
\end{equation*}
where we used Sobolev embeddings and
\eqref{eq:systhyp}. 
For the second term we use 
\begin{equation*} 
\( S \d_{t}  \d_{x}^{\alpha } \bu^\h , \d_{x}^{\alpha } \bu^\h\)
= \frac{\h}{2}\(S L( \d_{x}^{\alpha } \bu^\h )  , \d_{x}^{\alpha }
\bu^\h\)
  -  
\Big( S  \d_{x}^{\alpha } \Big(\sum _{j=1}^{n} A_j(\bu^\h) \d_{j}
  \bu^\h  \Big) , 
 \d_{x}^{\alpha } \bu^\h\Big).
\end{equation*} 
We notice that $SL$
is a skew-symmetric second order operator, so the first term is zero. 
The second term can be rewritten under the form
\begin{equation*} 
\begin{aligned}
\Big( S  \partial _{x}^{\alpha } \Big(\sum _{j=1}^{n} A_j(\bu^\h) \d_j
 \bu^\h  \Big)& ,   
 \partial _{x}^{\alpha } \bu^\h\Big)
= 
\Big(S  \sum _{j=1}^{n} A_j(\bu^\h) \d_j \partial _{x}^{\alpha
 } \bu^\h, 
 \partial _{x}^{\alpha } \bu^\h\Big) \\
 +& \Big(S \Bigl(  \partial _{x}^{\alpha } (\sum _{j=1}^n
 A_j(\bu^\h) \d_j \bu^\h  ) - 
  \sum _{j=1}^n A_j(\bu^\h) \d_j \partial _{x}^{\alpha }
 \bu^\h  \Bigr ), 
 \partial _{x}^{\alpha } \bu^\h\Big).
\end{aligned}
\end{equation*}
By symmetry of $SA_j(\bu^\h)$, 
\begin{equation*}
\begin{aligned}
\Big(S \sum _{j=1}^{n} A_j(\bu^\h) \d_j  \partial_{x}^{\alpha } \bu^\h, 
 \partial _{x}^{\alpha } \bu^\h\Big)
 = & - \sum _{j=1}^{n} \Big(\d_j (S A_j(\bu^\h))  \d_{x}^{\alpha } \bu^\h,
 \d_{x}^{\alpha } \bu^\h\Big)\\
& - \sum _{j=1}^{n} \Big(S A_j(\bu^\h) \d_j  \partial _{x}^{\alpha } \bu^\h,
 \partial _{x}^{\alpha } \bu^\h\Big)\, .
\end{aligned} 
\end{equation*} 
Therefore, so long as $\|\bu^\h\|_{L^\infty}\le 2\|a_0\|_{L^\infty}$
for $t\le 2$, 
\begin{equation*} 
\left| \Big(S  \sum _{j=1}^n A_j(\bu^\h) \d_j \partial _{x}^{\alpha } \bu^\h,
 \d_{x}^{\alpha } \bu^\h\Big) \right|
 \lesssim
 \left\|\d_{x}^{\alpha } \bu^\h \right\|_{L^{2}}^{2} 
\left\| \nabla _{x} \bu^\h \right\|_{L^{\infty }} \lesssim
 \left\|  \bu^\h \right\|_{H^s}^3 \, .
\end{equation*}
The usual estimates on commutators (see e.g. \cite{Majda}) lead to
\begin{equation*}
\left|\bigg(S  \Bigl (  \partial _{x}^{\alpha } \Big(\sum _{j=1}^n
  A_j(\bu^\h) \d_j \bu^\h  \Big) - 
  \sum _{j=1}^n A_j(\bu^\h) \d_j \partial _{x}^{\alpha } \bu^\h \Bigr ),
 \partial _{x}^{\alpha } \bu^\h\bigg) \right|
\le C\(\left\|\bu^\h\right\|_{H^s}\) \left\|\bu^\h\right\|_{H^s}^{2} .
\end{equation*}
Notice that $S^{-1}$ can be bounded by $C(\|\bu\|_{H^s})$, thus we
have proved:
\begin{equation*}
\frac{d}{dt} \sum _{|\alpha | \le s} \(S \partial _{x}^{\alpha } \bu^\h,
\partial _{x}^{\alpha } \bu^\h\) \le C\(\left\|\bu^\h\right\|_{H^s}\)
\sum _{|\alpha | \le s} \(S \partial _{x}^{\alpha } \bu^\h,
\partial _{x}^{\alpha } \bu^\h\)\, ,
\end{equation*}
for $s > 2+d/2 $.
Gronwall lemma along with a
continuity argument yield the counterpart of
\cite[Theorem~1.1]{Grenier98}: 
\begin{prop}\label{prop:p}
 Under Assumptions~\ref{hyp:tout}  with $n>k>0$, let
$s>2+n/2$. Then there exist $T>0$ independent of $\h\in ]0,1]$  and
$\psi^\h(t,x) = a^\h(t,x) e^{i\phi^\h(t,x)/\h}$ solution to \eqref{eq:psi} on
$[t_0^\h,T+t_0^\h]$. Moreover, $a^\h$ and
$\phi^\h$ are bounded 
in $L^\infty([t_0^\h,T+t_0^\h];H^s)$,
uniformly in $\h\in ]0,1]$. 
\end{prop}
\subsection{Convergence and small time properties}
\label{sec:end}
We can now complete the proof of Proposition~\ref{prop:bkwreduit}. For
$s>2+n/2$, we know that $\widetilde a^\h$ and $\widetilde {\tt v}^\h$
are bounded in 
$L^\infty([0,T];H^s)$, uniformly in $\h\in ]0,1]$. Using
\eqref{eq:systhyp}, we infer that $\d_t \widetilde a^\h$ and
$\d_t\widetilde {\tt v}^\h$ are 
bounded in $L^\infty([0,T];H^{s-2})$. Therefore, a subsequence of
$(\widetilde a^\h,\widetilde {\tt v}^\h)$ converges 
uniformly in  $C([0,T];H^{s'}_{\rm loc})$ to $(a,{\tt v})$ solution
of \eqref{eq:bkwlimite} for any $s'<s-1$. Decreasing $T$ if necessary,
we know that \eqref{eq:bkwlimite} has a unique solution in
$C([0,T];H^s)$ ( see e.g. \cite{Majda}).
By uniqueness for \eqref{eq:bkwlimite}, the whole sequence
$(\widetilde a^\h,\widetilde {\tt v}^\h)$ is convergent. 

For $t_0^\h\le t\le T$, write
\begin{equation*}
 {\tt v}^\h(t,x)= \widetilde {\tt v}^\h(t,x) -\int_{t-t_0^\h}^t \d_t
 \widetilde {\tt v}^\h(s,x)ds \ ;\  a^\h(t,x)= \widetilde a^\h(t,x)
 -\int_{t-t_0^\h}^t \d_t  
 \widetilde a^\h(s,x)ds 
\end{equation*}
Still using the boundedness of  $\d_t \widetilde a^\h$ and
$\d_t\widetilde {\tt v}^\h$ in $L^\infty([0,T];H^{s-2})$, we deduce
that the sequence  $(a^\h,
{\tt v}^\h)$ also converges to \eqref{eq:bkwlimite}. To estimate the
convergence and prove that it holds not only in $H^s_{\rm loc}$ , but
also in $H^s$, we shall use the analysis of
\eqref{eq:bkwlimite} for small times.
\begin{prop}\label{prop:small}
  Let Assumptions~\ref{hyp:tout} be satisfied, with $n>k>1$. Then
  there exist sequences $(\phi_j)_{j\ge 1}$ and $(a_j)_{j\ge 1}$ in
  $H^\infty(\R^n)$, such that the solution of
  \eqref{eq:bkwlimite} satisfies
  \begin{equation*}
    \phi(t,x)\sim \sum_{j\ge 1}t^{jn-1}\phi_j(x)\ ,\text{ and}\quad a(t,x)\sim
    \sum_{j\ge 0}t^{jn}a_j(x)\quad \text{as }t\to 0\, . 
  \end{equation*}
\end{prop}
Plugging such asymptotic series into \eqref{eq:bkwlimite}, a formal
computation yields a source term which is $\O(t^\infty)$ as $t\to
0$. The first terms are computed
 in \eqref{eq:phi1}--\eqref{eq:phi2}. Note that the series for $a$
 starts with $j=0$: the notations are consistent. 
The result then follows with the same approach as in the proof of
Proposition~\ref{prop:p},  
and  Borel lemma (see e.g. \cite{RauchUtah}). 
We can now be more precise about the convergence of $(a^\h,\phi^\h)$:
\begin{prop}\label{prop:estprec}
Let $s\in \N$. There exists $C$ independent of $\h$ such that for
every $t_0^\h\le t\le T$, 
\begin{equation*}
  \| a^\h (t)- a(t)\|_{H^s}
  +\| \phi^\h(t) - \phi(t)\|_{H^s} \le C\( \h t +
  \h^{\frac{\g(n-1)}{1-\g}}\). 
\end{equation*}
\end{prop}
\begin{proof}
  We keep the same notations as in the previous subsection,
  \eqref{eq:systhyp}. Define $(\widetilde a, \widetilde {\tt v})$ from
  $(a,{\tt v})$ by the same shift in time as above. Denote by $\bu$
  the analog of $\bu^\h$ 
  corresponding to $(\widetilde a,\widetilde {\tt v})$. We have
  \begin{equation*}
\d_t \(\bu^\h-\bu\) +\sum_{j=1}^n A_j(\bu^\h)\d_j \(\bu^\h-\bu\) +
\sum_{j=1}^n \(A_j(\bu^\h)-A_j(\bu)\)\d_j \bu  
  = \frac{\h}{2} L 
  \bu^\h    \, .
  \end{equation*}
Keeping the symmetrizer $S$ corresponding to $\bu^\h$, we can do
similar computations to those of the previous paragraph. Note that we
know that $\bu^\h$ and $\bu$ are bounded in $L^\infty([0,T];H^s)$. Denoting ${\bw}^\h
=\bu^\h-\bu$, we get, for $s>2+n/2$: 
\begin{equation*}
\frac{d}{dt} \sum _{|\alpha | \le s} \(S \partial _{x}^{\alpha } \bw^\h,
\partial _{x}^{\alpha } \bw^\h\) \lesssim
\sum _{|\alpha | \le s} \(S \partial _{x}^{\alpha } \bw^\h,
\partial _{x}^{\alpha } \bw^\h\)+ \h \| \bw^\h(t)\|_{H^s}\, .
\end{equation*}
Using Gronwall lemma and Proposition~\ref{prop:small}, we infer:
\begin{equation*}
  \| \bw^\h(t)\|_{H^s} \lesssim \h t  + \left\| (a,{\tt
      v})\big|_{t=t_0^\h}- (a,{\tt
      v})\big|_{t=0}\right\|_{H^s}\lesssim \h t + \(t_0^\h\)^{n-1}\,  .
\end{equation*}
This completes the proof of Proposition~\ref{prop:estprec}.
\end{proof}
So far, we have not used the assumption $k>1$. It appears when one
wants to approximate $e^{i\phi^\h/\h}$ by $e^{i\phi/\h}$: the factor
$1/\h$ requires some care. Proposition~\ref{prop:estprec} shows that 
for $t_0^\h\le \tau\le T$, 
\begin{equation*}
  \sup_{t_0^\h\le t\le \tau}\left\|
  \phi^\h(t)-\phi(t)\right\|_{H^s (\R^n)}= \O\(
  \h^{\frac{\g(n-1)}{1-\g}}\)+ \h \O(\tau)\, .   
\end{equation*}
Recalling 
that $\g = k/n$, we then have, for $s>n/2$,  
\begin{equation*}
\begin{aligned}
  \left\| \psi^\h - a_0e^{i\phi/\h}\right\|_{L^\infty([t_0^\h,\tau];H^s)}
  \lesssim &  
\left\| a^\h - a\right\|_{L^\infty([t_0^\h,T];H^s)} +
\left\| a - a_0\right\|_{L^\infty([t_0^\h,\tau];H^s)} \\
&+
\frac{1}{\h}\left\| \phi^\h -
  \phi\right\|_{L^\infty([t_0^\h,\tau];H^s)} \\
\le &\  \O(\h) +\O(\tau^n)+ \O\(\h^{\frac{k-1}{1-\g}}\) + \O(\tau)\, .
\end{aligned}
\end{equation*}
This completes the proof of Proposition~\ref{prop:bkwreduit},
and Theorem~\ref{theo:main} follows.  

\begin{rema}[Well-prepared data] If we had $t_0^\h=0$, then the
  assumption $k>1$ could be weakened to $k>0$. Back to the transform
  \eqref{eq:pseudo}, if we assume that
  \begin{equation*}
    u^\e_{\mid t=0} = a_0(x)e^{-i\frac{|x|^2}{2\e}} \exp\(i\e^{\g-1}
    \phi\(\e^\g,x\)\)\,,
  \end{equation*}
then the $\O\(\h^{\frac{k-1}{1-\g}}\)$ term in the above estimate
disappears, and we can conclude as before, supposing only
$n>k>0$ (but still $n\ge 2$). Recall that if $0<k\le 1$, 
then nonlinear effects are relevant at leading order for any positive
time (nonlinear propagation); they show up precisely in the phase
$\phi$.  
\end{rema}
\begin{rema}
  The cascade of phase shifts can be understood as the creation of a
  new phase, appearing discretely in time. With the transform \eqref{eq:pseudo}
  in mind, the asymptotic expansion of the phase shift $\phi$ stems
  from Proposition~\ref{prop:small}. The coupling in
  \eqref{eq:bkwlimite} shows that even if $\phi_{\mid t=0}=0$,
  $\phi(\delta,x)$ is not identically zero, for $\delta >0$
  arbitrarily small. The phase of  $\psi^\h$ is
  given asymptotically by
  \begin{equation*}
    \frac{\phi(t,x)}{\h} \sim \sum_{j\ge 1} \frac{t^{jn-1}}{\h}\phi_j(x)\,.
  \end{equation*}
  With the same line of reasoning as in the introduction, a phase
  shift appears for $t$ of order $\h^{1/(n-1)}$, then a second for $t$
    of order $\h^{1/(2n-1)}$, and so on. The superposition of these
      phase shifts, which are oscillating faster and faster, finally
      leads to a continuous phase, corresponding to an oscillation
      associated to the wavelength $\h$. The idea that the cascade of
      phase shifts corresponds to the (discrete in time) creation of a new
      (continuous) phase is reinforced by Equation~\eqref{eq:chg} below.
\end{rema}

\begin{rema}\label{rema:bkwameliore}
  Like in \cite{Grenier98}, we can prove a WKB type asymptotics, not
  only for small times as in Proposition~\ref{prop:bkwreduit}, but on
  the whole interval $[0,T]$. Indeed, with the notations of the proof
  of Proposition~\ref{prop:estprec}, $\bw^\h =
  \O(\h)$ in $L^\infty([0,T];H^s)$ for any $s\in \N$. Reasoning as
  above, $\bw^\h/\h $ is bounded in  $L^\infty([0,T];H^s)$,
  hence $\d_t \bw^\h/\h$ is bounded in
  $L^\infty([0,T];H^{s-2})$, and a subsequence converges to the
  linearization of \eqref{eq:systhyp} about $\bu$. By uniqueness, the
  whole sequence is convergent, and we have an error estimate in the
  same spirit as in Proposition~\ref{prop:estprec}. Denoting ${\bf d}
  = \lim \frac{\phi^\h-\phi}{\h}$, we have:
  \begin{equation*}
    \left\| \psi - a e^{i{\bf
    d}}e^{i\phi/\h}\right\|_{L^\infty([0,T];H^s)} \Tend \h 0 0 \, .
  \end{equation*}
This is a WKB type asymptotics, with the amplitude $a e^{i{\bf
    d}}$ and the phase $\phi$. 
\end{rema}

\section{Stability issues}
\label{sec:stab}

The construction of Section~\ref{sec:formal} and the results of the
previous paragraph do not agree. To see this, we come back to
Proposition~\ref{prop:small}: in \eqref{eq:bkwlimite}, we have
\begin{align}
\O\(t^{n-2}\):\quad &  \phi_1(x) = \frac{1}{n-1}f'(0)|a_0(x)|^2\,
,\label{eq:phi1} \\
\O\(t^{n-1}\):\quad & a_1 + \nabla \phi_1\cdot \nabla a_0
+\frac{1}{2}a_0\Delta \phi_1 =0\, ,\label{eq:a1}\\
\O\(t^{2n-2}\):\quad &  (2n-1)\phi_2 +\frac{1}{2}|\nabla \phi_1|^2
+2\operatorname{Re}(\overline{a_0}a_1)f'(0)
+\frac{f''(0)}{2}|a_0|^4=0\, . \label{eq:phi2}
\end{align}
The function $\phi_1$ is the same as the one obtained by the approach
of Section~\ref{sec:formal}: the two approximate solutions are close
to each other up to the first boundary layer, when the first phase
shift appears. On the other hand, we see that to get
$\phi_2$, the modulation of the amplitude ($a_1$) must be taken into
account; in \eqref{eq:evolg}, $g_2$ is computed without
evaluating $\Delta a_0$, unlike $\phi_2$. This means in particular
that the two approximate solutions diverge when reaching the second
boundary layer: the approach of Section~\ref{sec:formal} is only
formal, and does not lead to a good approximation. And yet, the source
term in Proposition~\ref{prop:formal} is small: thus, the linearized
semi-classical Schr\"odinger operator is not stable, in the
semi-classical limit. We will see below that this instability is not
due to a spectral instability, but to the fact that the approach
followed to construct the formal approximation was too crude.

This phenomenon is due to the super-criticality of the problem.
Indeed, for fixed $\e$, we 
deal with a nonlinear Schr\"odinger equation with repulsive
nonlinearity ($f'>0$), for which global well-posedness results are
available (see Remark~\ref{rema:reg}). When using the transform
\eqref{eq:pseudo}, notice that the 
parameter $\h$ in \eqref{eq:psi} goes to zero as $\e\to 0$ only when
$n>k$, that is in the super-critical case (compare with
Section~\ref{sec:critical}).

To understand better the instability mechanism, let us go back to the
comparison between the construction of 
Section~\ref{sec:formal} and the results of the previous
paragraph. Letting $N\to +\infty$ in Proposition~\ref{prop:formal}, we
have an approximate solution of the form
\begin{align*}
  v^\e (t,x) &= \frac{e^{i\frac{|x|^2}{2\e(t-1)}}}{(1-t)^{n/2}}
    a_0\(\frac{x}{1-t}\) \exp\(i\frac{1-t}{\e}
     g\(\frac{\e^\g}{1-t} \virgp \frac{x}{1-t}\)\)\\
&= \frac{e^{-i\frac{|\xi|^2}{2\h \tau}}}{(1-t)^{n/2}}
    a_0\(\xi\) \exp\(\frac{i}{\h }
     \frac{g\(\tau,\xi\)}{\tau}\)\Big|_{(\tau,\xi)=\(\frac{\e^\g}{1-t} \virgp
    \frac{x}{1-t}\) }.
\end{align*}
This formula and the transform \eqref{eq:pseudo} show that the
approximation of Section~\ref{sec:formal} is too crude, since it ignores
the coupling between phase and amplitude for \eqref{eq:psi}. 
Proposition~\ref{prop:small} and \eqref{eq:bkwlimite} show that to
have a good approximation of the phase, the coupling between phase and
amplitude must be taken into account at \emph{every} order. 

We can go one step further in the understanding of this apparent
instability, by applying the 
transform \eqref{eq:pseudo} to the intermediary approximate solution
$v_N^\e$. We show that the formal approximation stops being a good
approximation between the first and the second boundary layer. Assume
$\si=1$ so that the homogeneous nonlinearity satisfies
Assumptions~\ref{hyp:tout}. 
Like for the exact solution, write
\begin{equation*}
  v_N^\e(t,x) =
  \frac{1}{(1-t)^{n/2}}\psi_N^\e\(\frac{\e^\g}{1-t}\virgp
  \frac{x}{1-t} \) e^{i\frac{|x|^2}{2\e(t-1)}}\, .
\end{equation*}
Using the expression \eqref{eq:r_N}, we check that $\psi_N^\h$
solves
\begin{align*}
  i\h \d_t \psi_N^\h +\frac{\h^2}{2}\Delta \psi_N^\h = t^{n-2}
   |\psi_N^\h|^2 
  \psi_N^\h +\theta_N^\h(t,x)\, ,
\end{align*}
along with the initial condition
$\psi_N^\h\big|_{ t =
  \h^{\frac{\g}{1-\g}}} = a_0(x)+\O\(\h^{(\alpha -1)(1-\g)}\)$ in
$H^s(\R^n)$ for any $s>0$, where:
\begin{equation*}
  \theta_N^\h(t,x) = \(t^{(N+1)n-2}K_0(x) + i\h K_1(t,x)\)
  \psi_N^\h(t,x) +i\h K_2(t,x) + \h^2 K_3(t,x)\, ,
\end{equation*}
for some ``nice'' functions $K_j$. Now write $\psi_N^\h(t,x) =
a_N^\h(t,x) e^{i\phi_N^\h(t,x)/\h}$. We have:
\begin{equation}
\label{eq:systhyp2}
\begin{aligned}
  \d_t \bv^\h &+\sum_{j=1}^n A_j(\bv^\h)\d_j \bv^\h
  = \frac{\h}{2} L 
  \bv^\h + {\tt S}^\h (t,x)\, , \text{  with }\bv^\h(t,x) = \left(
    \begin{array}[l]{c}
      \operatorname{Re} a_N^\h \\
      \operatorname{Im} a_N^\h\\
      \d_1\phi_N^\h  \\
      \vdots \\
      \d_n\phi_N^\h
    \end{array}
\right),\\
  \text{and}&\quad {\tt S}^\h(t,x) = (t+t_0^\h)\left(
    \begin{array}[l]{c}
      K_1+ {\rm Re}\( (K_2-i\h K_3)e^{i\phi_N^\h/\h}\) \\
      K_1+ {\rm Im}\( (K_2-i\h K_3)e^{i\phi_N^\h/\h}\) \\
      -(t+t_0^\h)^{(N+1)n-3}\d_1K_0  \\
      \vdots \\
      -(t+t_0^\h)^{(N+1)n-3}\d_nK_0
    \end{array}
\right)\, ,
\end{aligned}
\end{equation}
where the matrices $A_j$ are the same as in Section~\ref{sec:grenier} and the
functions in the definitions of $\bv^\h$ and ${\tt S}^\h$ are evaluated
at $(t+t_0^\h,x)$. We can proceed like in Section~\ref{sec:grenier}:
the new term is the source ${\tt S}^\h$. Unlike for the exact
solution, the oscillatory aspect of the problem has not disappeared:
the first two components of ${\tt S}^\h$ contain a highly oscillatory
factor. Therefore, we cannot expect $\h$ independent energy
estimates here. To measure the effect of this oscillatory term, forget
the shift in time, and take $ t_0^\h=0$. Then assuming that for
small times, $\d_x^a \phi_N^\h (t,x) = \O(t^{n-1})$ for any
multi-index $a$ (like for the exact solution), the $H^s$ norms of
the first two components of ${\tt S}^\h$ are controlled by
\begin{equation*}
  \O\( t + \frac{t^{1+s(n-1)}}{\h^s}\).
\end{equation*}
A source of order $\O(t)$ is not a problem, since we eventually
consider the limit $t\to 0$. On the other hand, let us examine the
last term. Back to the initial variables, this yields a control by
\begin{equation*}
  \( \frac{\e^\g}{1-t}\)^{1+s(n-1)} \e^{-s(1-\g)} = \frac{\e^{\g
  +s\alpha -s}}{(1-t)^{1+s(n-1)}}\, .
\end{equation*}
This is small for $1-t\gg \e^\omega$, with
\begin{equation*}
  \omega = \frac{\g
  +s\alpha -s}{1+s(n-1)}\, .
\end{equation*}
We check that for $n >\alpha =k>1$, we have
\begin{equation*}
  \beta = \frac{\alpha -1}{n-1}< \omega = \frac{\g
  +s(\alpha -1)}{1+s(n-1)} < \frac{2\alpha -1}{2n-1}\, , \text{ for
  any } s\ge 1\,.
\end{equation*}
The first inequality means that we can expect the formal approximation
to be a good approximation of the exact solution beyond the first
boundary layer (and indeed, it is close to the approximate solution of
Section~\ref{sec:justif}). The second one explains why the
approximation ceases to 
be relevant before the second boundary layer. 

A possible way to understand the above computation is that the choice
of the variables is crucial: working with the ``usual'' unknown $v^\e$
(as in Section~\ref{sec:formal}) is not very efficient. On the other
hand, with 
the variables introduced by E.~Grenier for his generalized WKB
methods, a precise and rigorous 
analysis is possible, \emph{via} the transform
\eqref{eq:pseudo}. Thus, adding new variables helps the analysis:
this goes in the same direction as the general theory of geometric
optics, and the recent approach followed by
C.~Cheverry for a refinement of this principle 
\cite{CheverryCMP,CheverryBullSMF}.

\section{After the cascade of phase shifts}
\label{sec:beyond}

As announced in the introduction, our analysis stops for times of
order $t= 1-\lambda\e^\g$. For $\lambda \to +\infty$, we have
Theorem~\ref{theo:main}. For bounded $\lambda$, $\lambda \in
[\frac{1}{T},+\infty[$, the first part of
Proposition~\ref{prop:bkwreduit} provides an asymptotic
description; see also Remark~\ref{rema:bkwameliore}. This shows in
particular that the  
solution $\psi^\h$ is approximated in terms of  nonlinear geometric
optics: the eikonal equation contains the amplitude, therefore the
geometry of propagation needs not be the same as before, which occurred
along rays
$\frac{x}{1-t}=\text{const}$. Note also that the transform
\eqref{eq:pseudo} changes the space variable into a parameterization
of the family of rays, when they are straight lines. This explains
Figure~\ref{fig:rays}. Moreover, between the initial time $t=0$ and
$t = 1-\frac{\e^\g}{T}$, the order 
of magnitude of $u^\e$ changes.  Indeed,
Proposition~\ref{prop:bkwreduit} shows that for $t\in [t_0^\h,T]$,
$\psi^\h$ is of order $\O(1)$ in $L^\infty(\R^n)$ (take $s>n/2$ and
use Sobolev 
embeddings). By \eqref{eq:pseudo}, we infer that the amplitude of
$u^\e$ varies like $(1-t)^{-n/2}$, and changes from $\O(1)$ initially,
to $\O(\e^{-k/2})$ for $1-t\thickapprox \e^\g$. Such an amplification
is similar to what happens in the linear case.

The semi-classical conformal transform  \eqref{eq:pseudo}
cannot be interesting for values of $t$ too close to $1$, since it
becomes singular. It seems reasonable to introduce the ($L^2$ unitary)
scaling transform,
\begin{equation}
  \label{eq:transf2}
  u^\e(t,x)=\frac{1}{\e^{n\g/2}}\varphi^\e\( \frac{t-1}{\e^\g}\virgp
  \frac{x}{\e^\g}\)= \frac{1}{\e^{k/2}}\varphi^\e\( \frac{t-1}{\e^\g}\virgp
  \frac{x}{\e^\g}\) .
\end{equation}
With the same change of notation as for $\psi$ in
Section~\ref{sec:pseudo}, we have
\begin{equation}\label{eq:varphi}
  i\h \d_t \varphi^\h +\frac{\h^2}{2}\Delta \varphi^\h = f\( 
  |\varphi^\h|^2\) 
  \varphi^\h\, .
\end{equation}
We now have exactly the same equation as in \cite{Grenier98}. On the
other hand, let us examine the initial condition. Taking into account
the data $u^\e(0,x)=a_0(x)$, we find
\begin{equation*}
  \varphi^\h \( -\h^{\frac{\g}{\g -1}}, x\)=
  \h^{\frac{k}{2(1-\g)}}a_0\( \h^{\frac{\g}{1-\g}}x\) , 
\end{equation*}
that one may try to decouple to
\begin{equation*}
  \varphi^\h(t,x) \sim \frac{1}{|t|^{n/2}}a_0\(\frac{x}{t}\)
  e^{i\frac{|x|^2}{2\h t}}\quad \text{as }t\to -\infty \text{ and
  }\h\to 0\, .
\end{equation*}
Of course, the above limits $t\to -\infty$ and $\h\to 0$ do not
commute. Denote $U^\h(t)$ the 
unitary group associated to the linear semi-classical Schr\"odinger
equation ($f\equiv 0$ in the above equation). Then for fixed $\h>0$, 
\begin{equation*}
  U^\h(t)\varphi_0(x)\Eq t {-\infty} \frac{1}{(2i\pi \h
  t)^{n/2}}\widehat \varphi_0\(\frac{x}{\h t}\) e^{i\frac{|x|^2}{2\h
  t}}\, .
\end{equation*}
Fix $\h$ in the above asymptotics. If $f$ has moderate growth as
in Remark~\ref{rema:reg} (cubic nonlinearity in space dimension two or
three for instance), then there is scattering for
\eqref{eq:varphi} with $\h$ fixed, and
\begin{equation*}
  \varphi^\h(0,x)\Eq \h 0 \frac{1}{\h^{n/2}}\Phi\(\frac{x}{\h}\)\, ,
\end{equation*}
for some concentrating profile $\Phi$, where the powers of $\h$ stem
from the asymptotics for the free operator $U^\h$. We saw in the
introduction that 
at least when the nonlinearity is homogeneous, the conservation of
energy rules out such a possibility, in the limit $\h\to 0$.  

It is probably more interesting to try to match 
with the results of Proposition~\ref{prop:bkwreduit}. Comparing
\eqref{eq:pseudo} and \eqref{eq:transf2}, Proposition~\ref{prop:p}
yields, for
$\frac{-1}{t_0^\h}\le t\le\frac{-1}{T+t_0^\h}$:
\begin{equation}\label{eq:chg}
\begin{aligned}
  \varphi^\h(t,x) &= \(\frac{-1}{t}\)^{n/2}\psi^\h \(\frac{-1}{t}\virgp
  \frac{-x}{t} \) e^{i\frac{|x|^2}{2\h t}}\\
& =\(\frac{-1}{t}\)^{n/2}
  a^\h \(\frac{-1}{t}\virgp \frac{-x}{t} \) 
 \exp \( \frac{i}{\h}
\( \frac{|x|^2}{2 t} +\phi^\h\(\frac{-1}{t}\virgp \frac{-x}{t} \)
\) \)  \\
& =: {\bf a}^\h(t,x)e^{i{\Phi}^\h(t,x)/\h} \, .
\end{aligned}
\end{equation}
Note however that the term $\frac{|x|^2}{2 t}$ in the
phase does not belong to any Sobolev space; we would have to adapt the
statement of Proposition~\ref{prop:p} before claiming that we have
$\h$ independent estimates. We can then try to use Grenier's ideas
again to extend the lifespan of ${\bf a}^\h$
and $\Phi^\h$ to $[\frac{-1}{T+t_0^\h},{\widetilde T}]$, with
suitable $\h$ independent estimates.
We shall not pursue this point of view here. 

We conclude this section by
listing a series of questions that remain:

Do we have a WKB like description of $u^\e$ for some time $t>1$? If yes,
with one or several phases?

Is there a caustic for $\varphi^\h$? This is not even
  clear. Indeed, the initial problem \eqref{eq:princ} contains a data
  which causes focusing at one point in the linear case. However, we
  saw above that when reaching the boundary layer of size $\e^\g$,
  phase and amplitude of the solution become coupled in such a way
  that the geometry of the propagation is modified. If by any chance
  there is no caustic for $\varphi^\h$, then we might take $\widetilde T$
  arbitrarily large, and hope to get a description of $u^\e$ for any
  time. 

Note also that when WKB
  asymptotics is valid for $\varphi^\h$, then the nonlinear term in
  the conservation of the energy reaches its maximal order of
  magnitude. This is how we found the parameter $\g$ in the
  introduction for an homogeneous nonlinearity. In the case of
  Assumptions~\ref{hyp:tout}, things are similar. Denote
  \begin{equation*}
    F(y)=\int_0^y f\(\eta^2\)\eta d\eta\quad (\eta\text{ is a real
    variable}). 
  \end{equation*}
Then the generalization of the conservation of energy in
\eqref{eq:conservation} is: 
\begin{equation*}
  \frac{1}{2}\|\e\nabla_x u^\e(t)\|_{L^2}^2 + \e^{-k}\int_{x\in\R^n}
  F\(\e^{k/2} |u^\e(t,x)|\)dx  = \text{const.}=\O (1) \Eq \e 0
\frac{1}{2}\|x a_0\|_{L^2}^2\, .   
\end{equation*}
The nonlinear term in the above energy is exactly $\int F\(|\varphi^\h
\(\frac{t-1}{\e^\g},x\)|\)dx$, and is of order $\O(1)$ for,
say,
$\frac{t-1}{\e^\g} \in
[ \frac{-2}{T},\frac{-1}{T}]$, where WKB
  asymptotics for $\varphi^\h$ stems from
  Proposition~\ref{prop:bkwreduit}. (For $ \frac{t-1}{\e^\g} \to 
  -\infty$, dispersive properties of $\varphi^\h$ 
  make the nonlinear term small.)
On the other hand, we saw that if there is a caustic for
  $\varphi^\h$, then it cannot be reduced to a (single)
  point. 

Finally, the apparent instability discussed in
  Section~\ref{sec:stab} suggests that computing reliable numerical
  simulations to understand the asymptotic behavior 
  of $u^\e$ is a challenging problem, even before a caustic is formed,
  if there is any. Understanding the behavior of
  $(a^\h,\phi^\h)$ and $({\bf a}^\h,\Phi^\h)$, rather than working on
  $u^\e$ directly,  would certainly be more
  reasonable.

\bibliographystyle{amsplain}
\bibliography{../carles}

\end{document}